# ROBUST UNDERSTANDING OF WORD PROBLEMS WITH EXTRANEOUS INFORMATION

Yefim Bakman

Tel Aviv University, E-mail: bakman@post.tau.ac.il

## Abstract

Understanding of free-format multi-step arithmetic word problems with extraneous information is discussed. A model including a full set of general skills necessary for understanding such problems was developed and computer implemented. The validity of the simulation was confirmed by testing on a variety of word problems with extraneous information. The results of this study show the way for moving from an intuitive understanding of word problems to a conscious and controlled process.

## 1. Introduction

The main goal of the present study was elucidation of the full set of skills necessary for correct understanding of arithmetic word problems containing extraneous information. The problems were restricted to those whose solution required only two operations: addition and/or subtraction. Thus, the goal of the study was construction of a model of the most basic skills of understanding addressing no special knowledge.

Only few publications have been devoted to understanding word problems with extraneous information. Low and Over's (1989, 1990) experimental work showed that performance scores on different types of tasks (identifying irrelevant or missing information, problem solving, and problem categorization) were highly correlated. Muth (1992) investigated how the presence of extraneous information in word problems reduced the accuracy of children's solutions. Littlefield and Rieser's (1993) study is the only publication about strategies children use to cope with irrelevant information in a special kind of one-step word problems in which a superset is the requested quantity and the solver has to identify subsets of this superset in the problem text. This skill constitutes an integral part of the full understanding ability but it is not general enough to work in other cases.

On the other hand, in the domain of one-step word problems without extraneous information many articles have been published, which have developed two major theoretical models. The first model was developed in the studies of Riley, Greeno, and Heller (1983), and was further elaborated by Kintsch and Greeno (1985), Cummins, Kintsch, Reusser and Weimer (1988). The second model, called CHIPS, was suggested by Briars and Larkin (1984). The two models categorize simple word problems into the same three categories: compare, combine, and change (see Table 1). Both models were computer-implemented (see Dellarosa, 1986; Fletcher, 1985; Briars and Larkin, 1984) and their ability to solve one-step problems with only one possible operation - addition/subtraction - was proven.

The two simulations, Briars and Larkin's CHIPS, and Dellarosa's (1986) ARITHPRO, impose strict constraints on change problems: they can handle only one change verb, i.e. "to give". However there are many other change verbs denoting change of place (e.g. to put, to move), creation or termination of objects (e.g. to build, to eat), and problems containing them are unsolvable in these simulations. Another limitation of the two simulations concerns the change problem format: the first sentence of the problem must describe the number of objects the owner had to begin with, whereas the second sentence must contain the verb *gave*. Any variation of this format will make the problem unsolvable in the simulations.

The present work uses the earlier studies as a starting point. The analysis of the process of understanding presented in this article shows that the two change schemas offered by the previous studies do not supply an understanding of word problems that comprise change verbs other than "to give". At least eight different change schemas are needed to differentiate between situations involving changes in place, ownership, creation of new objects, and termination of existing objects. Some change verbs have complex meaning because they denote more than one change. Preliminary splitting of sentences comprising such verbs is necessary. Additional understanding skills were introduced in the model reported here so that the above restrictions were eliminated. Free-format multi-step word problems with extraneous information and containing any change verbs have become solvable in the present simulation.

Another goal of this study was to address the question of how skillful solvers manage to handle word problems with large amounts of irrelevant information. Experiments with adults and school students were conducted to test the hypothesis that proficient solvers employ a special strategy when deciding on schema relevancy (article in preparation). The experiments confirmed that good solvers consider a schema to be relevant if all its components have been found among the problem propositions. Using this strategy, proficient solvers managed to verbally solve word problems with extraneous information within two minutes whereas less skillful subjects were busy writing a system of equations that often led to errors.

The cautious strategy was implemented in the simulation ROBUST. The simulation parses natural language. It was successfully tested on a variety of multi-step word problems with extraneous information. The theoretical model of robust understanding was developed through experiments with the computer simulation of robust understanding in the process of development and through testing human subjects (students and adults). The resulting model and the experience of testing robust understanding might help in the development of instructional methods and in the diagnosis of robust understanding of word problems. This in turn enables students' intuitive understanding of word problems to evolve into a conscious and controlled process.

## 2. To understand understanding

According to Webster's Dictionary (1974) the word "understanding" has three meanings: a process, a mental state, and an ability. In order to differentiate between them we will call them *the process of understanding, the state of understanding,* and *the ability to understand. The process of understanding* is a process of selecting an adequate model for representing an object/phenomenon. Adequacy occurs if there exists an analogical correspondence between the model elements and those of the object/phenomenon. If such a match is successful, the person or the intelligent system that has made the selection is in *a state of understanding,* which means that he/she/it has represented the object/phenomenon by the model.

It follows that *the ability to understand* in a certain domain of knowledge includes a set of models and procedures supplying representations for objects/phenomena of the domain. It is clear that the quality of understanding depends on the variety of the models and on the perfection of the procedures for model selection. The human ability to understand develops and is improved over time, with every individual going through a novice stage when his/her ability to understand is immature. For example, younger children may use the following simple method for word problem solving: "If the word *altogether* occurs in the problem then add the two quantities given in the problem". In spite of its primitiveness the method contains the main blocks of understanding: the procedure providing model selection (the search for the key word *altogether*), and the model itself, comprising two given quantities and one unknown, which equals their sum. Even though it is primitive the method results in the correct solutions for many simple arithmetic word problems containing the word *altogether.*

In our considerations we will use the term *schema* to denote an adjustable mental structure containing variables (Rumelhart, 1980). If constants are substituted for the variables in the structure, then we obtain a model which constitutes an instantiation of the schema. Thus each schema helps to generate an infinite number of models.

According to Riley et al's (1983) hypothesis, four schemas are enough to model children's understanding of simple word problems. This approach was further developed by Kintsch and Greeno (1985), and Cummins et al. (1988), who called the four schemas *Change-In, Change-Out, Compare*, and *Combine*. Each of the schemas includes three quantities (variables) interrelated by the equality $a + b = c$. A single instantiation of one of the schemas is sufficient to represent a one-step arithmetic word problem. Hence, such problems can be classified into three categories: change problems, compare problems and combine problems (see Table 1). Since we are looking at multi-step word problems, a list of instantiations of the schemas is needed to represent such problems.

**Table 1. Examples of the three problem types.**

| Problem types | Example problems |
|---|---|
| *Change* | John had 5 apples. Mary gave him 3 apples. How many apples does John have now? (Change-In) |
| | Mary had 5 plums. Then she gave 3 plums to Tom. How many plums does Mary have now? (Change-Out) |
| *Combine* | Ruth has 3 dolls. Ann has 4 dolls. How many dolls do they have altogether ? |
| *Compare* | Sara has 6 flowers. Clara has 3 flowers more than Sara. How many flowers does Clara have ? |

### 3. Robust understanding

Transpositions of problem sentences and/or addition of irrelevant data to the problem text must not affect the problem solution. If the ability of a solver to understand word problems is stable to such changes in the problem text then we will call this kind of understanding **robust understanding**. A teacher who tests his/her students by well-specified word problems will not be able to discriminate between immature and robust understanding. On the other hand, word problems with extraneous information discriminate between the two levels of understanding. To illustrate this let us consider the following problem with extraneous information:

> *Problem 1*. Ruth had 3 apples. She put 2 apples into a basket. How many apples are there in the basket now, if in the beginning there were 4 apples in the basket?

To solve change problems Dellarosa (1986) placed covert constraints on the possible problem format: the first sentence must be of the type "Owner A had $m$ objects". The second sentence format must be either "A gave B $n$ objects" or "B gave A $n$ objects". The choice of the change schema is made according to the rule "If the actor of the second sentence is A, then a Change-Out schema is triggered, otherwise, a Change-In schema is brought into play". Applied to Problem 1, this rule yields that a Change-Out schema fits the problem because the owner in the first sentence is Ruth and she is the actor in the second sentence. But in fact, the actor Ruth is irrelevant in solving the problem. How can the solver find out this fact?

To make the choice between the Change-In and Change-Out schemas Briars and Larkin (1984) suggested the following method: if the change verb is followed by "him" or "her", this means that the owner has more objects now. Otherwise he/she has less. Again we see that the method does not work in the case of Problem 1. Littlefield and Rieser's (1993) model searches through the problem text for the actor, action, and unit of measure that match those requested in the question. This strategy can be applied only to a special type of word problems.

Our approach says that a proficient solver must distinguish between change verbs which denote change of ownership (e.g., *to obtain, to lose*) and between those denoting change of place (e.g., *to put, to move*). More examples of such verbs are presented in Table 2. Knowing that the verb *to put* refers to a change of place, a skillful solver will ignore the actor Ruth of Problem 1 and will concentrate on the basket (the place). The latter will help the solver to decide which of the two sentences "Ruth had 3 apples" or "There were 4 apples in the basket" is relevant to the change expressed in the phrase "Ruth put 2 apples in the basket".

**Table 2. Change verb categorization.**

| Types of changes | Change verbs |
|---|---|
| Transfer-In-Ownership | to receive, to get |
| Transfer-Out-Ownership | to lose, to forfeit, to send |
| Transfer-In-Place | to fetch, to bring, to put in, to lay, to enter, to fall into, to add |
| Transfer-Out-Place | to take out, to take away, to exit, to go away, to send, to drag out, to fall from |
| Creation | to build, to be born, to create, to make |
| Termination | to eat, to destroy, to die, to kill |

Note. The verbs *to buy, to give, to pay, to sell, to donate, to steal* include two simple meanings each because they involve two persons, one of which gets something while the other forfeits the same.

A confirmation of this can be found in Schank and Abelson's (1977) study, which discussed types of primitive acts, among which are ATRANS (transfer of an abstract relation such as possession, ownership or control) and PTRANS (transfer of physical location). Schank and Abelson also indicated that some verbs (e.g., *to buy*) are made up of more than one primitive act. The same idea is presented in the note to Table 2 where compound change verbs are listed. To this group belongs the most common change verb *to give,* which denotes two elementary changes at once, "to forfeit" and "to get".

Change verbs denoting transfer of ownership/location do not exhaust all actions which result in changes of amounts. There are two other groups of such verbs. Let us consider the following phrase: "Tom built two houses in the village". Though the two houses have not been transferred into the village, it is clear that the number of the village houses has increased by two. There must be an owner or owners of the houses but we cannot be sure that this is Tom. Thus, there are verbs (e.g., to build, to be born) whose meaning is creation of new objects in a certain place and, possibly, of a certain ownership so that a change of the corresponding amounts occurs. Similarly, there is another category of change verbs (e.g., to break, to eat) denoting the inverse process, namely, "termination" of objects, which results in diminishing the amount of the corresponding objects in a certain place and, perhaps, in someone's possession. The full categorization of change verbs is presented in Table 2. According to this categorization we have at least six change situations: Transfer-In-Ownership, Transfer-Out-Ownership, Transfer-In-Place, Transfer-Out-Place, Creation, Termination.

The previous discussion shows that unavailability of such a differentiation will prevent the solver to build adequate representations of word problems with extraneous information. Since the earlier simulations do not categorize change verbs, the skills they present do not supply a full understanding of word problems under consideration. The present article is meant to make up for this deficiency.

### 4. Formulas as a means for recognition of change schema instantiations

As was shown in the previous section, the categorization of change verbs helps to determine whether a change situation refers to change in place and/or to change of ownership, but the initial and final amounts of change can be scattered through the problem text. The earlier simulations avoided this difficulty by stipulating that the initial and the final amounts must be recorded in the first and the third sentences of the problem respectively. But here we are discussing free-format word problems, and therefore we need a special procedure to search for, and to collect, the three quantities of a change schema.

To accomplish this we will use a new notion, i.e. "change formula", which constitutes a generalized description of a change situation. There is exactly one change formula for each change schema. For example, the following formula corresponds to the Transfer-In-Place schema (other change formulas are presented in Table 3):

> There was an initial number of objects in a place.
> An additional number of objects was transferred into this place.
> There is a final number of objects in the place.

The words "objects", "initial number", "place" and so on are variables. Their names are not important (e.g., we may write X instead of "initial number") because concrete values will be assigned to the variables later on. The idea is to draw a correspondence between members of a sentence containing a change verb, and the formula variables. For example, if we have the sentence "2 apples were put into the basket", then the correspondence is as follows:

> objects = apples
> place = basket
> additional number = 2

Substitution of the obtained values of the variables into the two other phrases of the Transfer-In-Place formula yields the following <u>formula instantiation</u>:

> There was an initial number of apples in the basket.
> 2 apples were transferred into the basket.
> There is a final number of apples in the basket.

Thus, the instantiation supplied us with two new sentences which describe the two other amounts of the change schema. Those two sentences obtained are "There was an initial number of apples in the basket" and "There is a final number of apples in the basket". Now it is easy to find the descriptions of the initial and final amounts of the schema among the problem sentences. If we deal with Problem 1, they are "There were 4 apples in the basket" (i.e., initial number = 4) and "How many apples are there in the basket" which yields the final number = ? (i.e., the required quantity).

The corresponding interrelation between the three quantities constitutes the Transfer-In-Place schema instantiation which we will record in the following concise form: Transfer-In-Place (initially 4, in 2, finally ?).

The equality between the three quantities is known. It is ? = 4 + 2.

**Table 3. Formulas used by ROBUST for change schema recognition**

| Name of change schema | Formula associated with the schema |
|---|---|
| Transfer-In-Place | There were *X* objects in the place. *Y* objects were transferred into the place. There are *Z* objects in the place now. |
| Transfer-Out-Place | There were *X* objects in the place. *Y* objects were transferred out of the place. There are *Z* objects in the place now. |
| Transfer-In-Ownership | The owner had *R* objects. The owner received *S* objects. The owner has *T* objects now. |
| Transfer-Out-Ownership | The owner had *R* objects. The owner forfeited *S* objects. The owner has *T* objects now. |
| Creation (ownership) | The owner had *R* objects. The owner created *S* objects. The owner has *T* objects now. |
| Creation (place) | There were *X* objects in the place. *Y* objects were created in the place. There are *Z* objects in the place now. |
| Termination (ownership) | The owner had *R* objects. The owner terminated *S* objects. The owner has *T* objects now. |
| Termination (place) | There were *X* objects in the place. *Y* objects were terminated in the place. There are *Z* objects in the place now. |

## 5. The process of understanding

For our specific case of word problems, a state of understanding constitutes a possible representation of the problem by means of sets with schema instantiations interrelating between the sets. For example, Problem 1 can be represented by means of four sets expressed in the form of propositions (for the sake of simplicity we will write them as simple English sentences):

1) Ruth had 3 apples.

2) Ruth put 2 apples into the basket.

3) There are ? apples in the basket now.

4) There were 4 apples in the basket.

and one schema instantiation

Transfer-In-Place (initially 4, in 2, finally ?).

Now we shall recapitulate the main steps of the process of understanding as it was implemented in the ROBUST simulation. Suppose that parsing has been completed successfully, i.e. all the problem sentences have been turned into propositions. Alongside the parsing, ROBUST creates *the initial list of schema instantiations* which contains all instantiations of Compare and Combine schemas found in the problem text. Since this step is similar to its analogue in Dellarosa's (1986) simulation, we do not describe it in detail. Then the procedure of understanding is as follows:

1. Search among the problem propositions for those having complex sense and split them into elementary ones.
2. For each elementary change verb and the proposition containing the verb:
   a. Choose the single change formula related to the chosen change verb and draw a correspondence between the formula variables and the matching constant values in the proposition.
   b. Substitute the constant values for the correspondent formula variables in all propositions of the formula to produce a formula instantiation.
   c. There are two formula propositions which have not been matched yet to those of the problem. Search for their instantiations among the propositions of the problem.
   d. **If** the search was successful **then**

      substitute the newly matched constant values for the correspondent variables in all propositions of the formula instantiation to update the latter. Create the Change schema instantiation corresponding to the formula instantiation. Attach the former to the *list of schema instantiations (LSI)* representing the problem.

      **end-if**

**end-of-procedure**

Now we will apply the procedure to the following sample problem:

   *Problem 2.* David gave 3 candies to Ruth, and John gave 2 candies to David. Now David has 4 candies more than Ruth has. How many candies does David have now, if Ruth had 7 candies in the beginning ?

First of all the problem text undergoes parsing to turn all the sentences into the following propositions (for details see section 8):

   1) David gave 3 candies to Ruth.
   2) John gave 2 candies to David.
   3) David has 4 candies more than Ruth has.
   4) David has ? candies. (the question)

5) Ruth had 7 candies.

Proposition 3 constitutes a comparison, hence it gives rise to the instantiation of Compare schema: MORE (?, than X, by 4), where the new unknown X is defined by the additional proposition:

6) Ruth has X candies.

The initial list of schema instantiations (LSI) contains one instantiation: MORE (?, than X, by 4). ROBUST's output at this point is shown in Table 4.

**Table 4. ROBUST's output after parsing Problem 2**

| Propositions | Schema Instantiations |
|---|---|
| David gave 3 candies to Ruth | More (?, than X, by 4) |
| John gave 2 candies to David | |
| David has ? candies | |
| Ruth had 7 candies | |
| Ruth has X candies | |

Now step 1 begins, according to which each complex change verb must be split into its elementary components. There are two sentences with the complex verb "to give" in our example. The elementary components are "to forfeit" and "to get". So we receive the following set of propositions:

1a) David forfeited 3 candies.

1b) Ruth got 3 candies.

2a) John forfeited 2 candies.

2b) David got 2 candies.

3) David has 4 candies more than Ruth has.

4) David has ? candies (the question).

5) Ruth had 7 candies.

6) Ruth has X candies.

Step 2. Each change verb in the propositions is examined as to whether or not the correspondent change schema is relevant to the problem. The first proposition "David forfeited 3 candies" comprises the change verb "forfeited" which belongs to the category Transfer-Out-Ownership. Its formula is:

An owner had an initial number of objects. An additional number of objects was transferred from the owner. The owner has a final number of objects.

To draw a correspondence between the second sentence of the formula and the proposition "David forfeited 3 candies" the latter should be presented in the canonical form in which the direct object will become a subject:

"3 candies were transferred from David".

Now the correspondence is easy to set. It yields

the owner = David

the objects = candies

the additional number = 3

Substituting the found values into the formula (step 2b) gives the following instantiation of the formula:

David had the initial number of candies.

3 candies were transferred from David.

David has the final number of candies.

Now step 2c begins. There are two formula propositions which have not yet been matched to those of the problem. We find that the formula proposition "David has the final number of candies" can be matched to proposition 4 "David has ? candies" but no problem proposition fits the first proposition of the formula instantiation. Since not all amounts of the schema are present among the problem propositions, a schema instantiation is not attached to LSI (the reason is explained in the next section). The same happens with propositions 2a and 2b. Only proposition 1b gives rise to a full Transfer-In-Ownership formula instantiation:

Ruth had 7 candies.

3 candies were transferred to Ruth.

Ruth has X candies.

Now that all formula propositions are matched to their instantiations in the problem text, we can create the corresponding schema instantiation:

Transfer-In-Ownership (initially 7, in 3, finally X).

This schema instantiation is attached to LSI (step 2d). Thus, only one of the four candidate change schemas appears to be relevant. The process of understanding results in the state of understanding - the problem is represented by means of two schema instantiations:

MORE (?, than X, by 4)

Transfer-In-Ownership (initially 7, in 3, finally X).

The procedure is not simple, and not every student can be expected to be able to develop it without assistance. Perhaps a lack of this skill is what prevents less successful students from succeeding in understanding word problems. But the skill is not necessary for solving simple word problems when two amounts are given and the task is to choose one of two possible operations (+ or -).

## 6. Cautious strategy for schema instantiation creation

If not all amounts of a change schema have been found in the problem, the solver has to decide whether to create a schema instantiation by introducing new unknowns or not. The *total* strategy is a strategy which creates a schema instantiation and attaches it to LSI every time an elementary change verb is encountered. This strategy may seem expedient at a first glance. Indeed, the person who uses it does not need to search through the problem text for two other number sets of the change schema, thus saving on effort. On the other hand, such a strategy fosters a superabundance of schema instantiations representing the problem, especially if the problem includes irrelevant information. The extra schema instantiations become a hindrance to problem understanding because they exhaust short-term memory.

According to the present model of robust understanding, proficient problem solvers use a kind of *cautious* strategy for creation of schema instantiations that warns not to record a schema instantiation if it is unclear whether this schema is necessary for the solution or not. Of course, schema relevancy can not be fully ascertained until the problem has been solved, and therefore the solver has to use a heuristic method to estimate schema relevancy at the earlier stage.

Simulation ROBUST uses the following version of the cautious strategy: a change schema is considered relevant (and its instantiation is recorded) if all its amounts are present in the problem text. Taking the example of Problem 2, the procedure of understanding described in the previous section creates no change schemas associated with David's candies. The reason is that not all quantities of the schemas are present among the problem propositions. And indeed, those change schemas are irrelevant for the solution.

If a problem includes extra information, use of the cautious strategy is decisive as to whether the solver can handle this problem or not. For such word problems a strong correlation between the strategy used and solution correctness was observed in experiments with human subjects (article in preparation).

The algorithm of problem solving will be published in another article.

**Table 5. Some of the problems solved by ROBUST (translated from Hebrew).**

---

1. Tom and Ruth had 8 apples altogether. Ruth gave Tom 3 apples. Now Tom has 5 apples. How many apples did Ruth have in the beginning?

2. Two boys left a room. 3 girls and 5 boys remained in the room. How many boys were there in the room in the beginning?

3. 5 girls bought 6 tickets. 7 boys bought 8 tickets. How many tickets did the children buy altogether?

4. Ruth had 5 nuts more than Dan had. Ruth gave Dan 3 nuts. Dan gave 2 nuts to David. Now Dan has 4 nuts and David has 6 nuts. How many nuts does Ruth have now?

5. In the beginning there were 4 eggs more in a refrigerator than there were in a box. Sara transferred 5 eggs from a basket into the refrigerator. After this David transferred 6 eggs from the basket into the box. 3 eggs fell out of the box. Now there are 2 eggs more in the box than there are in the basket. How many eggs are there in the basket now, if there are 12 eggs in the refrigerator?

6. Fred had 10 candies. Dan gave 2 candies to Susan and Fred gave 3 candies to Dan. Now Dan has 4 candies less than Susan has. How many candies does Dan have now, if it is known that Susan had 7 candies in the beginning and Fred has 9 candies now?[1]

---

## 7. Contribution

The difficulty of mapping a word problem on a change schema consists in the dispersal of the schema amounts through the problem text. In addition, the initial and the final amounts of a change schema are described differently for different change verbs. Two earlier simulations of arithmetic word problem solving CHIPS (Briars and Larkin, 1984) and ARITHPRO (Dellarosa, 1986) used immature strategies that could not cope with the indicated difficulties. In order for the simulations to handle word problems, rigid limitations were imposed on the order of the problem sentences. Denise Cummins pointed out that odd formats (e.g., Jill gave Mary 2 marbles. Jill had 2 marbles before that) perplexed children as well as ARITHPRO. Nor was it possible to diversify problems by substituting the verb "to give" by some other change verbs, say "to put".

---

[1] ROBUST detected a contradiction in the problem data.

The present model removes these limitations by using change formulas to search for and to collect the three quantities of a change schema. Now the order of the problem sentences may be arbitrary and any change verb may occur in the text.

Another important contribution of the present work is its implication for the category of compound change verbs. The verb "to give" also belongs to this category - it signals two possible schema instantiations at once. This possibility was not taken into account by the earlier models.

The simulation does not record all possible schema instantiations but rather eliminates irrelevant ones, thus copying the cautious strategy used by proficient problem solvers to cope with irrelevant data. Multi-step word problems can be solved as a consequence of the achieved robust understanding since quantities irrelevant to any particular step do not confuse the simulation.

The last contribution of the simulation is that it parses natural language.

## ACKNOWLEDGMENTS

This study was supported by the Science & Technology Educational Center of Tel-Aviv University and by the Science Absorption Center, the Ministry of Absorption, Israel. I am grateful to David Chen, Dina Tirosh, and David Miodusar for their assistance and discussions during the work, and to Denise Cummins and an anonymous reviewer for their helpful comments.